\documentclass{amsart}
\usepackage{amssymb}
\usepackage{amsfonts,amsmath,amssymb,mathrsfs,verbatim}
\usepackage[dvips]{graphicx}
\usepackage{graphics}
\usepackage{psfrag}
\begin{document}
\newtheorem{thm1}{Theorem}[section]
\newtheorem{lem1}[thm1]{Lemma}
\newtheorem{rem1}[thm1]{Remark}
\newtheorem{def1}[thm1]{Definition}
\newtheorem{cor1}[thm1]{Corollary}
\newtheorem{defn1}[thm1]{Definition}
\newtheorem{prop1}[thm1]{Proposition}
\newtheorem{ex1}[thm1]{Example}
\newtheorem{alg1}[thm1]{Algorithm}


\title[toric ideals of graphs]{Graver degrees are not polynomially bounded by true circuit degrees}
\author{Christos Tatakis}
\address{Mitilini, P.O. Box 13, Mitilini (Lesvos) 81100, Greece}
\email{chtataki@cc.uoi.gr}
\author{Apostolos Thoma }
\address{Department of Mathematics, University of Ioannina,
Ioannina 45110, Greece }
\email{athoma@uoi.gr}
\thanks{}


\date{}

\dedicatory{}

\begin{abstract}
\par
Let $I_A$ be a toric ideal.
We prove that the degrees of the elements of the Graver basis of $I_A$
are not polynomially bounded by the true degrees of the circuits of $I_A$.
\end{abstract}

\maketitle


\section{Introduction}

\par  Let $A=\{\textbf{a}_1,\ldots,\textbf{a}_m\}\subseteq \mathbb{N}^n$
be a vector configuration in $\mathbb{Q}^n$ and
$\mathbb{N}A:=\{l_1\textbf{a}_1+\cdots+l_m\textbf{a}_m \ | \ l_i \in
\mathbb{N}\}$ the corresponding affine semigroup.  We grade the
polynomial ring $\mathbb{K}[x_1,\ldots,x_m]$ over an arbitrary field $\mathbb{K}$ by the
semigroup $\mathbb{N}A$ setting $\deg_{A}(x_i)=\textbf{a}_i$ for
$i=1,\ldots,m$. For $\textbf{u}=(u_1,\ldots,u_m) \in \mathbb{N}^m$,
we define the $A$-\emph{degree} of the monomial $\textbf{x}^{\textbf{u}}:=x_1^{u_1} \cdots x_m^{u_m}$
to be \[ u_1\textbf{a}_1+\cdots+u_m\textbf{a}_m \in \mathbb{N}A.\]
We denoted by $\deg_{A}(\textbf{x}^{\textbf{u}})$, while the usual degree $u_1+\cdots +u_m$ of
$\textbf{x}^{\textbf{u}}$ we denoted by $\deg(\textbf{x}^{\textbf{u}})$. The
\emph{ toric ideal} $I_{A}$ associated to $A$ is the prime ideal
generated by all the binomials $\textbf{x}^{\textbf{u}}- \textbf{x}^{\textbf{v}}$
such that $\deg_{A}(\textbf{x}^{\textbf{u}})=\deg_{A}(\textbf{x}^{\textbf{v}})$, see \cite{St}. For such binomials,
we set $\deg_A(\textbf{x}^{\textbf{u}}- \textbf{x}^{\textbf{v}}):=\deg_{A}(\textbf{x}^{\textbf{u}})$.
A nonzero binomial $\textbf{x}^{\textbf{u}}- \textbf{x}^{\textbf{v}}$ in $I_A$ is called \emph{ primitive} if there exists no other binomial
 $\textbf{x}^{\textbf{w}}- \textbf{x}^{\textbf{z}}$ in $I_A$
such that $\textbf{x}^{\textbf{w}}$ divides
$ \textbf{x}^{\textbf{u}}$ and $\textbf{x}^{\textbf{z}}$
divides $ \textbf{x}^{\textbf{v}}$.
 The set of the primitive binomials forms the Graver basis of $I_A$ and is denoted by
$Gr_A$.
An irreducible binomial is called a \emph{ circuit} if it has minimal support. The set of the circuits is
denoted by $\mathcal{ C}_A$ and it is a subset of the Graver basis, see \cite{St}. One of the fundamental problems in toric
algebra is to give good upper bounds on the degrees of the elements of the Graver basis, see \cite{H, St, St1}.
It was conjectured that  the degree of any element in the Graver basis $Gr_A$ of a toric ideal $I_A$
is bounded above by the maximal true degree of any circuit in $\mathcal{ C}_A$, \cite[Conjecture 4.8]{St1},
\cite[Conjecture 2.2.10]{H}.
Following \cite{St1} we define the true degree of a
circuit as follows:
Consider any circuit $C\in \mathcal{ C}_A$ and regard its support supp($C$) as a  subset of $A$. The lattice $\mathbb{Z}$(supp($C$))
has finite index
in the lattice $\mathbb{R}$(supp($C$))$\cap \mathbb{Z}A$, which is called the index of the circuit $C$ and denoted by index($C$).
The \emph{ true degree} of the circuit $C$
is the product $\deg(C)\cdot $index($C$).
The crucial role of the true  circuit degrees was first highlighted in Hosten's dissertation \cite{H}.

Let us call $t_A$ the maximal true degree of any circuit in $\mathcal { C}_A$. The true circuit conjecture says
that $$\deg(B)\leq t_A,$$ for every $B\in Gr_A$.
There are several examples of families of toric ideals where the true circuit
conjecture is true, see for example
\cite{Pe}. The true circuit conjecture is also true for some families of toric ideals of graphs,
see \cite[Section 4]{TT}. However the true circuit conjecture is not true in the general case. In \cite{TT} we gave an infinite family of counterexamples
 to the true circuit conjecture by providing toric ideals and elements of the Graver basis
 for which their degrees are not bounded above by $t_A$.
 We note that in the counterexamples of \cite{TT} the degrees of the elements of the
 Graver basis were bounded by $ t_A^2$. In this article we consider the following question:

\textbf{Question:} {\em Does the degree of any element in the Graver basis $Gr_A$ of a toric ideal $I_A$
is bounded above by a constant times $(t_A)^2$ or a constant times $(t_A)^{2014}$?}

To disprove such a statement, one needs to compute the Graver basis and the set of circuits
for  toric ideals $I_A$ in a polynomial ring with a huge number of variables.
In order to produce  examples of toric ideals
such that there exist elements  in their Graver basis
 of very high degree and at the same time the true degrees of their circuits  have to be relatively low.
This procedure is computationally demanding,
if not impossible.
An alternative approach is given by the class of the toric ideals of graphs where
we explicitly know the form of the elements of their Graver basis, see \cite{ RTT}, and of their circuits,
see  \cite{Vi}.

The main result of the article is Theorem \ref{main} which says that

{\em there is no polynomial in $t_A$
that bounds the degree of any element in the Graver basis $Gr_A$ of a toric ideal $I_A$}.

To prove the theorem we are going to construct a family of examples of graphs $G_r^n$. For the toric ideals of these graphs and
for a fixed $n$
 we are going to prove that there are elements in the Graver basis
 whose degrees are exponential on $r$, see Proposition \ref{degprim}, while
 the true degrees of their circuits are linear on $r$, see Theorem \ref{truedeg} and Proposition \ref{degcircuit}.

\section{Toric ideals of graphs}\label{section 2}

\par
Let $G$ be a finite simple connected graph with vertices
$V(G)=\{v_{1},\ldots,v_{n}\}$ and edges $E(G)=\{e_{1},\ldots,e_{m}\}$.
Let $\mathbb{K}[e_{1},\ldots,e_{m}]$
be the polynomial ring in the $m$ variables $e_{1},\ldots,e_{m}$ over a field $\mathbb{K}$.  We
will associate each edge $e=\{v_{i},v_{j}\}\in E(G)$ with the element
$a_{e}=v_{i}+v_{j}$ in the free abelian group $ \mathbb{Z}^n $
with basis the set of vertices
of $G$.  Each vertex $v_j\in V(G)$ is associated with the vector $(0,\ldots,0,1,0,\ldots,0)$, where the
non zero component is in the $j$ position. We denote by $I_G$ the toric ideal $I_{A_{G}}$ in
$\mathbb{K}[e_{1},\ldots,e_{m}]$, where  $A_{G}=\{a_{e}\ | \ e\in E(G)\}\subset \mathbb{Z}^n $.

A \emph{walk}  connecting $v_{i_{1}}\in V(G)$ and
$v_{i_{s+1}}\in V(G)$ is a finite sequence of the form
$$w=(\{v_{i_1},v_{i_2}\},\{v_{i_2},v_{i_3}\},\ldots,\{v_{i_s},v_{i_{s+1}}\})$$
with each $e_{i_j}=\{v_{i_j},v_{i_{j+1}}\}\in E(G)$. A trail is a walk in which all edges are distinct.
The \emph{length}
of the walk $w$ is  the number $s$ of its edges. An
even (respectively odd) walk is a walk of \emph{even} (respectively odd) length.
A walk
$w=(\{v_{i_1},v_{i_2}\},\{v_{i_2},v_{i_3}\},\ldots,\{v_{i_s},v_{i_{s+1}}\})$
is called \emph{closed} if $v_{i_{s+1}}=v_{i_1}$. A \emph{cycle}
is a closed walk
$$(\{v_{i_1},v_{i_2}\},\{v_{i_2},v_{i_3}\},\ldots,\{v_{i_s},v_{i_{1}}\})$$ with
$v_{i_k}\neq v_{i_j},$ for every $ 1\leq k < j \leq s$.

Given an even closed walk $w$ of the graph $G$; where $$w =(e_{i_1}, e_{i_2},\dots,
e_{i_{2q}}),$$ we define
$$E^+(w)=\prod _{k=1}^{q} e_{i_{2k-1}},\ E^-(w)=\prod _{k=1}^{q} e_{i_{2k}}$$
and we denote by $B_w$ the binomial
$$B_w=\prod _{k=1}^{q} e_{i_{2k-1}}-\prod _{k=1}^{q} e_{i_{2k}}.$$
It is easy to see that $B_w\in I_G$. Moreover, it is known that the toric ideal $I_G$
is generated by binomials of this form, see \cite{Vi}. Note that the binomials $B_w$
are homogeneous and the degree of $B_w$ is $q$, the half of the number of edges of the walk.
For convenience,
we denote by $\textbf{w}$ the subgraph of $G$ with vertices the vertices of the
walk and edges the edges of the walk $w$. We call a walk
$w'=(e_{j_{1}},\dots,e_{j_{t}})$ a \emph{subwalk} of $w$ if
$e_{j_1}\cdots e_{j_t}| e_{i_1}\cdots e_{i_{2q}}.$
An even
closed walk $w$ is said to
be primitive if there exists no even closed subwalk $\xi$ of $w$ of smaller
length such that $E^+(\xi)| E^+(w)$ and $E^-(\xi)| E^-(w)$. The walk $w$
is primitive if and only if the binomial $B_w$ is primitive.

A \emph{ cut edge} (respectively \emph{ cut vertex}) is an edge (respectively vertex) of
the graph whose removal increases the number of connected
components of the remaining subgraph.  A graph is called \emph{
biconnected} if it is connected and does not contain a cut
vertex. A \emph{ block} is a maximal biconnected subgraph of a given
graph $G$.

The following theorems determine the form of the circuits and the
primitive binomials of a toric ideal of a graph $G$. R.~Villarreal in
\cite[Proposition 4.2]{Vi} gave a necessary and sufficient
characterization of circuits:
\begin{thm1}\label{circuit} Let $G$ be a graph and let $W$ be a connected subgraph of $G$.
The subgraph $W$ is the graph  ${\bf w}$ of a walk $w$ such that  $B_w$ is a circuit
 if and only if
\begin{enumerate}
  \item $W$ is an even cycle or
  \item $W$ consists of two odd cycles intersecting in exactly one vertex or
  \item $W$ consists of two vertex-disjoint odd cycles joined by a path.
\end{enumerate}
\end{thm1}

Primitive walks were first studied by T. Hibi and H. Ohsugi, see \cite{OH}. The next Theorem by  E.~Reyes, Ch.~Tatakis and A.~Thoma \cite{RTT} describes  the form of the underlying graph
of a primitive walk.

\begin{thm1} \label{primitive-graph}
Let $G$ be a  graph and let $W$ be a connected subgraph of $G$.
The subgraph $W$ is the graph  ${\bf w}$ of a primitive walk $w$
 if and only if \begin{enumerate}
  \item  $W$ is an even cycle or
  \item  $W$ is not biconnected and
\begin{enumerate}
  \item every block of $W$ is a cycle or a cut edge and
  \item every cut vertex of $W$ belongs to exactly two blocks and separates the graph in two parts, the total number of edges
of the cyclic blocks in each part is odd.
\end{enumerate}
\end{enumerate}
\end{thm1}

Observe that if $W'$ is the graph taken from $W$ by replacing every cut edge with two edges,
then $W'$ is an Eulerian graph since it is connected,  every cut vertex has degree  four and the others have degree two.
An  Eulerian trail   is a trail in a graph which visits every edge of the graph exactly once.
Any closed Eulerian trail
$w'$ of $W'$ gives rise to an even closed walk $w$ of $W$ for which every single edge
of the graph  $W'$ is a single edge of the walk $w$ and every multiple edge of the graph $W'$
is a double edge of the walk $w$ and a cut edge of $W={\bf w}$.
Different closed Eulerian trails may give different walks,
but all the corresponding binomials $B_w$ are equal or opposite.

\section{On the True circuit degree of  toric ideals of graphs}\label{Section 3}

In the next Theorem we prove that the index of any circuit $C$ in the toric ideal of a graph $G$ is equal to
1 and therefore the true degree of a circuit $C$ is equal to its degree.

\begin{thm1}\label{truedeg} Let $G$ be a graph and let $C$ be a circuit in $\mathcal{C}_{A_G}$.
Then $$true \deg(C)=\deg(C).$$
\end{thm1}
\textbf{Proof.} By definition $true \deg\ (C)=\deg(C)\cdot$ index($C$). We will prove that the  index($C$) is equal to one for every
circuit $C$ in a toric ideal of a graph $I_G$.
It is enough to prove that $\mathbb{Z}$(supp($C$))=$\mathbb{R}$(supp($C$))$\cap \mathbb{Z}A_G$.
Obviously $\mathbb{Z}$(supp($C$))$\subseteq \mathbb{R}$(supp($C$))$\cap \mathbb{Z}A_G$.
For the converse consider a circuit $C$ in $\mathcal{C}_{A_G}$.
By Theorem \ref{circuit} there are two cases.

First case: $C=B_w$  where ${\bf w}$ is an even cycle and let it be $$C=(e_1=\{v_{2k},v_1\},e_2=\{v_1,v_2\},\ldots,e_{2k}=\{v_{2k-1},v_{2k}\}).$$
Therefore supp($C$)$=\{ a_{e_1}, a_{e_2}, a_{e_3}, \ldots, a_{e_{2k}}\}$.
Since $C$ is a cycle we know that $$a_{e_1}-a_{e_2}+a_{e_3}-\ldots-a_{e_{2k}}=0. $$
Let ${\bf x}\in \mathbb{R}$(supp($C$))$\cap \mathbb{Z}A_G$, where $A_G=\{a_e|e\in E(G)\}$. Therefore ${\bf x}=r_1a_{e_1}+\ldots+r_{2k}a_{e_{2k}}$, where
$r_1,\ldots,r_{2k}\in\mathbb{R}$, and also  ${\bf x}\in \mathbb{Z}A_G\subset \mathbb{Z}^n$.
 By  ${\bf x}_{v}$ we denote the $v$ coordinate of  ${\bf x}$ in $\mathbb{Z}^n$ with the
 canonical basis denoted by the vertices of $G$.
 Then ${\bf x}_{v_1}=r_1+r_2\in \mathbb{Z}$, ${\bf x}_{v_2}=r_2+r_3\in \mathbb{Z}, \ldots, {\bf x}_{v_{2k}}=r_{2k}+r_1\in \mathbb{Z}$.
 It follows that $$r_{2l}\equiv-r_1\mod\mathbb{Z},\ \ r_{2l-1}\equiv r_1\mod\mathbb{Z},$$ for $1\leq l\leq k$.
 Therefore there exist integers $z_1=0, z_2,\dots ,z_{2k}$ such that
 $r_{2l}=z_{2l}-r_1$ and $r_{2l-1}=z_{2l-1}+r_1$. Then ${\bf x}=r_1a_{e_1}+\ldots+r_{2k}a_{e_{2k}}=
 r_1a_{e_1}+(z_2a_{e_2}-r_1a_{e_2})+(z_3a_{e_3}+r_1a_{e_3})+\ldots+(z_{2k}a_{e_{2k}}-r_1a_{e_{2k}})=
            z_2a_{e_2}+\ldots+z_{2k}a_{e_{2k}}\in \mathbb{Z}$(supp($C$)).

Second case: $C=B_w$  where ${\bf w}$  consists of two vertex disjoint odd cycles joined by a path or two odd cycles intersecting
in exactly one vertex, see Theorem \ref{circuit}.
Let $(e_1=\{v_{1},v_2\},e_2=\{v_2,v_3\},\ldots,e_{2l+1}=\{v_{2l+1},v_{1}\})$ be the one odd cycle, let
$(\xi_1=\{v_1, w_1\}, \xi_2=\{w_1, w_2\}, \ldots, \xi_{t}=\{w_{t-1}, u_1\}
)$ be the path of length $t$ and
$(\varepsilon_1=\{u_{1},u_2\},\varepsilon_2=\{u_2,u_3\},\ldots,\varepsilon_{2s+1}=\{u_{2s+1},u_{1}\})$
the second odd cycle. In the case that the length $t$ of the path is zero, $v_1=u_1$.
Therefore supp($C$)$=\{ a_{e_1}, a_{e_2},  \ldots, a_{e_{2l+1}}, a_{\xi_1}, \ldots,
a_{\xi_t}, a_{\varepsilon_1}, a_{\varepsilon_2},  \ldots,
a_{\varepsilon_{2s+1}}  \}$.
Since $C$ is a circuit we have that
$$ a_{e_1}-a_{e_2}  \ldots +a_{e_{2l+1}}-2a_{\xi_1}+ \ldots +2(-1)^ta_{\xi_t}+(-1)^{t+1} (a_{\varepsilon_1}-a_{\varepsilon_2}+ \ldots +
a_{\varepsilon_{2s+1}})=0. $$
Let ${\bf x}\in \mathbb{R}$(supp($C$))$\cap \mathbb{Z}A_G$ then
${\bf x}=r_1a_{e_1}+\ldots+r_{2l+1}a_{e_{2l+1}}+q_1a_{\xi_1}+ \ldots +q_ta_{\xi_t}+ \varrho_1a_{\varepsilon_1}+ \varrho_2a_{\varepsilon_2}+
\ldots+
\varrho_{2s+1}a_{\varepsilon_{2s+1}}  $, where
$r_1,\ldots,r_{2l+1}, q_1,\ldots, q_t, \varrho_1,\ldots, \varrho_{2s+1}\in\mathbb{R}$,
and also  ${\bf x}\in \mathbb{Z}A_G\subset \mathbb{Z}^n$.
 By looking at the  coordinates of  ${\bf x}$
 it follows that $$r_{2i}\equiv-r_1\mod\mathbb{Z},\ \ r_{2i+1}\equiv r_1\mod\mathbb{Z},$$
 $$q_{m}\equiv(-1)^m2r_1\mod\mathbb{Z},$$
 $$\varrho_{2j}\equiv(-1)^tr_1\mod\mathbb{Z},\ \ \varrho_{2j+1}\equiv (-1)^{t+1}r_1\mod\mathbb{Z},$$
for $1\leq i\leq l$, $1\leq m\leq t$ and $1\leq j\leq s$.
 Therefore there exist integers $x_2,\dots , x_{2l+1}, z_1,
\dots , z_t, w_1, \dots ,w_{2s+1}$ such that $r_j=x_j+(-1)^{j+1}r_1$, $q_j=z_j+2(-1)^{t+j}r_1$ and $\varrho_j=w_j+
(-1)^{t+j}r_1.$
  Then ${\bf x}=r_1a_{e_1}+\ldots+r_{2l+1}a_{e_{2l+1}}+q_1a_{\xi_1}+ \ldots +q_ta_{\xi_t}+ \varrho_1a_{\varepsilon_1}+ \varrho_2a_{\varepsilon_2}+
\ldots+
\varrho_{2s+1}a_{\varepsilon_{2s+1}}  =x_2a_{e_2}+\ldots+x_{2l+1}a_{e_{2l+1}}+z_1a_{\xi_1}+ \ldots +z_ta_{\xi_t}+
w_1a_{\varepsilon_1}+ w_2a_{\varepsilon_2}+
\ldots+
w_{2s+1}a_{\varepsilon_{2s+1}}  \in \mathbb{Z}$(supp($C$)). \\
Therefore in all cases $\mathbb{R}$(supp($C$))$\cap \mathbb{Z}A_G\subset \mathbb{Z}$(supp($C$)) and thus index$(C)=1$ for all
circuits $C$ in $I_{A_G}$. \hfill $\square$

\section{Bounds of Graver and True Circuit degrees}

The aim of this section is to provide examples of toric ideals
such that there are elements in their Graver bases that have very high
degree while the true degrees of their circuits remain relatively low.
We will do this for toric ideals of certain graphs, since the full power of
Theorem
\ref{truedeg} will come to use, and true degrees are equal to usual degrees.\\
Let $G_1$, $G_2$ be two vertex disjoint graphs,
on the vertices sets $V(G_1)=\{v_1,\ldots,v_s\}$, $V(G_2)=\{u_1,\ldots,u_k\}$
and on the edges sets $E(G_1), E(G_2)$ correspondingly.
We define the {\em sum of the graphs} $G_1,G_2$ on the vertices $v_i,u_j$ as a new graph $G$
formed from their union by identifying the pair of vertices  $v_i,u_j$ to form a single vertex $u$.
The new vertex $u$ is a cut vertex  in the new graph $G$ if both  $G_1$, $G_2$ are not trivial.
We say that we {\em add} to a vertex $v$ of a graph $G_1$ a cycle $S$, to get a graph $G$ if $G$ is the sum of
$G_1, S$ on the vertices $v\in V(G_1)$ and any vertex $u\in S$.

\begin{figure}
\begin{center}
\includegraphics[scale=2.0]{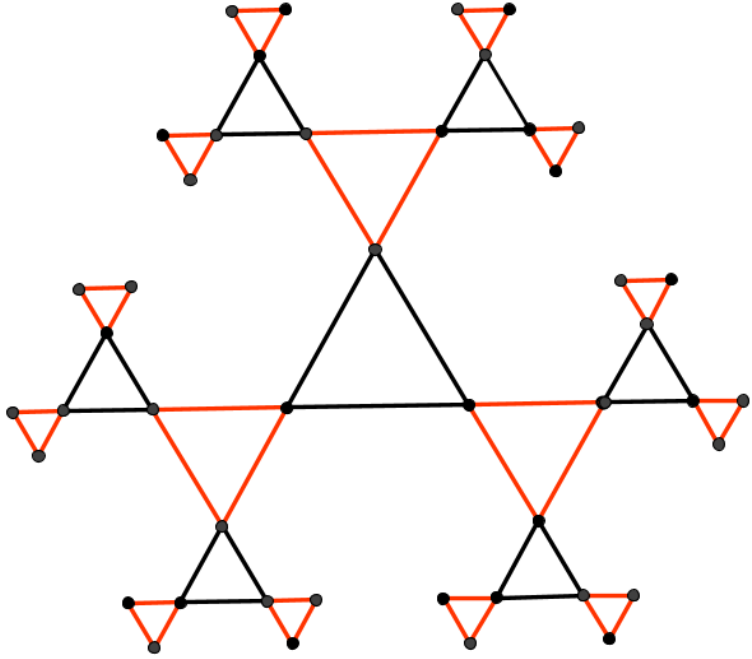}
\end{center}
\caption{The  graph $G_3^3$} 
\label{Figure 1}
\end{figure}

Let $n$ be an odd integer greater than or equal to three. Let $G_0^n$ be a cycle of length $n$.
For $r\ge 0$ we define the graph $G_{r}^n$ inductively on $r$.
$G_{r}^n$ is the graph taken from  $G_{r-1}^n$ by adding to each vertex of
degree two of the graph $G_{r-1}^n$  a cycle of length $n$.
Figure \ref{Figure 1} shows the  graph $G_3^3$.

We consider the graphs $G_0^n$ up to $G_{r-1}^n$ as subgraphs of $G_r^n$.
We note that the graph $G_r^n$ is Eulerian since by construction
it is connected and every vertex has even degree, four if it is also a vertex of $G_{r-1}^n$ and two if it  is not.

In the next Proposition we prove that the binomial $B_{w_r^n}$ belongs to the Graver basis of $I_{G_r^n}$ and compute its degree.

\begin{prop1} \label{degprim} Let $w_r^n$ be any closed Eulerian trail of the  graph $G_r^n$.
The binomial $B_{w_r^n}$ is an element of the Graver basis of $I_{G_r^n}$ and
$$\deg(B_{w_r^n})=\frac{1}{2}(n+n^2(\frac{(n-1)^r-1}{n-2})). $$
\end{prop1}

\textbf{Proof.}
We will prove the theorem by induction. We claim that the binomial $B_{w_s^n}$ belongs to the Graver basis of $I_{G_r^n}$, has degree $\frac{n+n^2(\frac{(n-1)^s-1}{n-2})}{2}$
and the graph $G_{s}^n={\bf w}_s^n$ has $n(n-1)^s$ vertices of degree 2, for
$1\leq s\leq r$.

For $s=1$ we consider the subgraph $G_1^n={\bf w}_1^n$ of $G_r^n$.
The graph is not biconnected, every block of the graph is a cycle and
there are no cut edges. Also every cut vertex of $G_1^n$ belongs to exactly two blocks
and separates the graph in two parts.
 One of them is a cycle of length $n$ and the other
 consists of $n$ cyclic blocks of $n^2$ total number of edges. Thus the total number of edges
of the cyclic blocks in each of the two parts is odd. Theorem \ref{primitive-graph} implies
 that $B_{w_1^n}$ is primitive. The total number of edges of $G_1^n$ is $n^2+n$, therefore
  the degree of the binomial $B_{w_1^n}$ is $\frac{n+n^2}{2}$ and the graph $G_{1}^n={\bf w}_1^n$ has $n(n-1)$ vertices of degree 2.
\\
Suppose that $B_{w_s^n}$ is primitive, $\deg(B_{w_s^n})=\frac{n+n^2(\frac{(n-1)^s-1}{n-2})}{2}$ and
the graph $G_{s}^n={\bf w}_s^n$ has $n(n-1)^s$ vertices of degree 2.
By the construction of the graph $G_{s+1}^n$, in every vertex of degree two of the graph $G_s^n$ we add an odd cycle of length $n$.
Since there are $n(n-1)^s$ vertices of degree two in $G_{s}^n$, the graph $G_{s+1}^n$  has $ n(n-1)^s$ new cycles, $n(n-1)^{s+1}$ vertices of degree 2
and $n\cdot n(n-1)^s$ new edges.
Therefore the binomial $B_{w_{s+1}^n}$ has degree
 \begin{eqnarray}
       \deg(B_{w_{s+1}^n}) &=&  \frac{n+n^2(\frac{(n-1)^s-1}{n-2})}{2}+\frac{n^2(n-1)^s}{2}\nonumber\\
                    &=&
                    \frac{n+n^2(\frac{(n-1)^{s+1}-1}{n-2})}{2}.\nonumber
\end{eqnarray}

The graph  $G_{s+1}^n={\bf w}_{s+1}^n$ is not biconnected
and every block of the graph is a cycle,
since the graph $G_{s+1}^n$ is constructed by adding cycles on the vertices of degree
two of the  graph $G_{s}^n$. Let $v$ be a cut vertex of the graph $G_{s+1}^n$.
The vertex $v$ is also a vertex of the subgraph $G_s^n$. There are two cases.
Either the vertex $v$ is a cut vertex of the subgraph $G_s^n$ or it has degree two in $G_s^n$.\\
First case, the vertex $v$ is a cut vertex in the graph $G_s^n$.
By the hypothesis $B_{w_s^n}$ is primitive,
therefore the vertex $v$ separates the graph $G_s^n={\bf w}_{s}^n$ in two parts.
The total number of edges
of the cyclic blocks in each of the two parts is odd by Theorem \ref{primitive-graph}.
 The graph $G_{s+1}^n$ is taken from the graph  $G_s^n$ by adding in every vertex of degree two of
 $G_s^n$ a cycle of length $n$. Thus
 in each cycle  of the graph
$G_s^n$ that has $n-1$ vertices of degree two we add $(n-1)n$ new edges,
i.e. even number of edges and therefore the vertex $v$
 separates also the graph $G_{s+1}^n$  in two parts, the total number of edges
of the cyclic blocks in each part is odd. \\
In the second case, the vertex $v$ has degree two in the graph $G_s^n$.
The vertex $v$  separates the graph $G_{s+1}^n$ in two parts.
One of them is a cycle of length $n$
and the other one has $2\deg(B_{w_{s+1}^n})-n$ edges. Thus the total number of edges
of the cyclic blocks in each part is odd.\\
From Theorem \ref{primitive-graph} we conclude that the binomial $B_{w_{s+1}^n}$ is primitive. \hfill $\square$

Let $B(G_r^n)$ be the {\em block tree} of $G_r^n$, the bipartite graph with bipartition
$(\mathbb{B},\mathbb{S})$ where $\mathbb{B}$ is the set of blocks of
$G_r^n$ and $\mathbb{S}$ is the set of cut vertices of $G_r^n$,
$\{{\mathcal B}, v\}$ is an edge if and only if $v\in {\mathcal B}$. The leaves of the block tree are the vertices of the block tree which have degree one.
Let ${\mathcal B}_k, {\mathcal B}_i, {\mathcal B}_l$ be blocks of a graph $G_r^n$. We call the block ${\mathcal B}_i$ \emph{internal block}
 of ${\mathcal B}_k,{\mathcal B}_l$, if ${\mathcal B}_i$ is an internal vertex in the unique path defined by ${\mathcal B}_k,{\mathcal B}_l$ in the block tree $B(G_r^n)$.
Every path of the graph $G_r^n$ from the block ${\mathcal B}_k$ to the block ${\mathcal B}_l$ passes
from every internal  block of ${\mathcal B}_k,{\mathcal B}_l$.
The path has vertices at least the cut vertices  which are vertices in
 the  path $({\mathcal B}_k,\ldots,{\mathcal B}_l)$
  in $B(G_r^n)$ and from one to at most $n-1$ common edges with the cycle that forms an internal block.

We denote by $\d({\mathcal B}_1,{\mathcal B}_2)$ the block distance between two vertices ${\mathcal B}_1,
{\mathcal B}_2\in \mathbb{B}$
of the block tree $B(G_r^n)$, which we define as the number of the internal vertices belonging to $\mathbb{B}$
in the unique path defined by the blocks ${\mathcal B}_1,{\mathcal B}_2$ in the block tree  $B(G_r^n)$.

The next lemma will be used to prove proposition \ref{degcircuit}.

\begin{lem1}\label{leaves} Let ${\mathcal B}_1,{\mathcal B}_2$ be two blocks of the graph $G_r^n$. Then $$\d({\mathcal B}_1,{\mathcal B}_2)\leq 2r-1.$$
\end{lem1}
\textbf{Proof.} We will prove it by induction. We claim that for any two blocks ${\mathcal B}_1,{\mathcal B}_2$
of the graph $G_s^n$ holds  $\d({\mathcal B}_1,{\mathcal B}_2)\leq 2s-1,$ for $1\leq s\leq r$. \\
We consider the block tree $B(G_1^n)$.  Let ${\mathcal B}_1,{\mathcal B}_2$ be two blocks of the graph $G_1^n$.
If both of them are leaves of the block tree $B(G_1^n)$ then $\d({\mathcal B}_1,{\mathcal B}_2)=1$ since there
is exactly one internal block, which corresponds to the graph $G_0^n$.
Otherwise, the distance is equal to 0. In every case $\d({\mathcal B}_1,{\mathcal B}_2)\leq 1=2\cdot1-1.$ \\
Suppose that the claim is true for $G_s^n$.
We consider the graph $G_{s+1}^n$ and let ${\mathcal B}_1,{\mathcal B}_2$ be two of its
blocks. Each of the blocks  ${\mathcal B}_1,{\mathcal B}_2$ is either block of the graph $G_s^n$
or has a common cut vertex with a block of the graph  $G_s^n$. It follows from
the induction hypothesis that $\d({\mathcal B}_1,{\mathcal B}_2)\leq (2r-1)+2=2(r+1)-1$.
 \hfill $\square$

We denote by $t_{A_{G_r^n}}$ the maximum degree of a circuit in the graph $G_r^n$.
In the following proposition we are providing a bound for the $t_{A_{G_r^n}}$.

\begin{prop1}\label{degcircuit} Let $t_{A_{G_r^n}}$ the maximum degree of a circuit in the graph $G_r^n$. Then $t_{A_{G_r^n}}\leq n+(2r-1)(n-1)$.
\end{prop1}
\textbf{Proof.} The graph $G_r^n$ has no even cycles and therefore the subgraph corresponding to a circuit
consists by two different  odd cycles joined by a path, see Theorem \ref{circuit}.  We remark that
every cycle of the graph has length $n$ and it is a block. Therefore it is enough to prove that a path between two
blocks ${\mathcal B}_1,{\mathcal B}_2$ of $G_r^n$ has length at most $(2r-1)(n-1)$.
Each such path passes from all internal blocks of ${\mathcal B}_1,{\mathcal B}_2$ and no other and has at most $n-1$ common edges  with
every one of them.
Therefore the path has at most length $\d({\mathcal B}_1,{\mathcal B}_2)\cdot (n-1)\leq (2r-1)(n-1).$ Thus the corresponding circuit has
degree at most $ n+(2r-1)(n-1)$. \hfill $\square$

\begin{rem1}{\rm It is not difficult to see that the bound given  at Proposition \ref{degcircuit}
 is sharp. In fact, there are several appropriate choices for the two blocks ${\mathcal B}_1,{\mathcal B}_2$ of $G_r^n$ and a unique
 choice of the path between them such that the $t_{A_{G_r^n}}=n+(2r-1)(n-1)$.}
\end{rem1}

 There are several bounds
on the degrees of the elements of the Graver basis of a toric ideal, see for example \cite{H, St, TT}.
The following theorem is the main result of the paper. It shows that for a general toric ideal $I_A$
 a bound given by a polynomial     in $t_A$ for the degrees of the elements of the Graver basis
 does not exist. Recall that $t_A$ is the
maximal true degree of a circuit in $I_A$.

\begin{thm1} \label{main} The degrees of the elements in the Graver basis  of a toric ideal $I_A$ cannot
be bounded polynomially above by the maximal true degree of a circuit.
\end{thm1}
\textbf{Proof.} Let $G$ be the graph $G_r^n$. It follows from Theorem \ref{truedeg} and Proposition \ref{degcircuit} that the maximal true
degree of a circuit is linear on $r$, while from Proposition \ref{degprim} there exists an element in the Graver basis
whose degree  is exponential in $r$. Therefore the degree of an element in the Graver basis $Gr_{A_G}$ of a toric ideal $I_{A_G}$ cannot
be bounded polynomially above by the maximal true degree of a circuit in $\mathcal{ C}_{A_G}$.
The proof of the theorem follows. \hfill $\square$

\end{document}